# Spider networks

**Leo Egghe[(1)], Li Li [(2)], Ronald Rousseau[(3),(4)]**

[(1)]Hasselt University, Hasselt, Belgium

E-mail: leo.egghe@uhasselt.be

ORCID: 0000-0001-8419-2932

[(2)] National Science Library, Chinese Academy of Sciences

33 Beisihuan Xilu, Zhongguancun, Beijing, P.R.China

E-mail: lili2020@mail.las.ac.cn

ORCID: 0000-0001-8326-3620

[(3)]KU Leuven, MSI, Facultair Onderzoekscentrum ECOOM,

Naamsestraat 61, 3000 Leuven, Belgium

E-mail: ronald.rousseau@kuleuven.be

&

[(4)]University of Antwerp, Faculty of Social Sciences,

Middelheimlaan 1, 2020 Antwerp, Belgium

E-mail: ronald.rousseau@uantwerpen.be

ORCID: 0000-0002-3252-2538

## Abstract

In this investigation we study a family of networks, called spiders, which covers a range of networks going from chains to complete graphs. These spiders are characterized by three parameters: the number of nodes in the core, the number of legs at each core node, and the length of these legs. Keeping two of the three parameters constant we investigate if spiders are small worlds in the sense recently defined by Egghe.

## 1. Introduction

### 1.1 Network notations and representations

Let G = (*V*,*E*) be an undirected network (or graph), where $V = (v_k)_{k=1,\dots,N}$ denotes the set of nodes or vertices and *E* denotes the set of links or edges. Networks can be represented in different ways, among which we will use the two-dimensional lay-out and adjacency matrices.

a) Two-dimensional graphs. By this, we mean a two-dimensional figure consisting of dots, representing nodes, and lines connecting the dots, representing the links. When two lines intersect, this has as such no meaning, in the sense that the intersection point is not necessarily a new node.

b) Adjacency matrices. These are square matrices for which the number of rows (hence also the number of columns) is equal to the number of nodes. The value of the cell (i,j) is 1 if node i is directly connected to node j, and 0 if this is not the case. The adjacency matrix of an undirected network is always symmetric.

### 1.2 Distance in a network

A path in a network is a sequence of different nodes one by one connected by edges. The distance $d_G$ (or simply d when it is clear in which network we work) between two nodes in the network G, is equal to the number of links situated on a shortest path, called a geodesic, between these two nodes. It is well-known that $d_G$ is a distance or metric in the mathematical sense of the word. Consequently, the distance between two nodes connected by an edge is equal to one. Each network studied in this article is assumed to be connected, i.e., there is a path between any two nodes. The maximum number of links in an N-node network is N(N-1)/2.

We recall for further use the definition of the following three basic networks. The chain, the star and the extended star, see Fig.1.

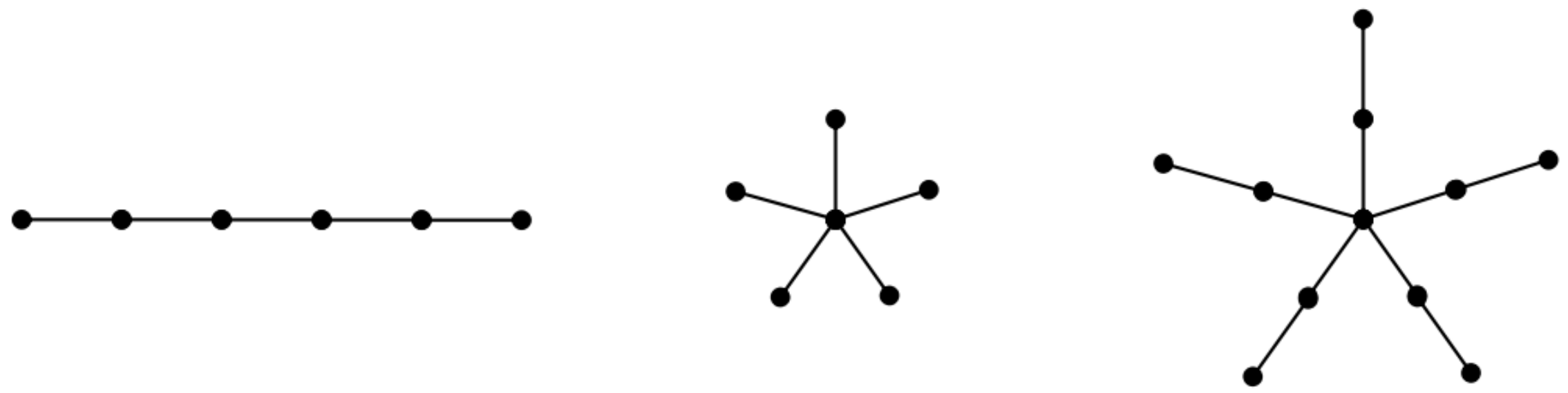

Fig.1 A chain, a star and an extended star

A chain is a graph G consisting of alternating vertices and edges, beginning and ending with vertices and in which each edge is incident with the two vertices immediately preceding and following it. A chain will be described through its number of edges. Hence an N chain has N+1 nodes. This description corresponds with the idea of six degrees of separation, where there are 6 links connecting 7 actors. A star is a graph consisting of a central node and some, say K, terminal nodes, connected only to the central node. Finally, an extended star has a fixed number of nodes between each terminal node and the central node. Hence, when there are K terminal nodes, and S nodes in between, there are (S+1)K edges, and in total N = (S+1)K+1 nodes. Fig. 1 shows the case S=1 and K=5.

In this study we only use undirected, unweighted, and connected networks.

## 2. Some network indicators

In this section, we recall some basic (classical) network indicators (Wasserman & Faust, 1994; Otte & Rousseau, 2002). These will be used further in this study.

The most basic network indicator is the density D(G). D(G) is defined as the number of links present, divided by the maximum number of links in a network with the same number of nodes. A network in which each node is connected to each other node is called a complete network (complete graph). The number of links in a complete N-node network is $\binom{N}{2} = \frac{N(N-1)}{2}$. Hence, if #E > 0 (and hence N > 1) denotes the number of links (edges) in the network G, then

$$D(G) = \frac{2\,(\#E)}{N(N-1)}\,. \tag{1}$$

Density is a property of the whole network. The next indicator refers to the role of a specific node in the whole network. The degree centrality of a node is defined as its number of neighbors. In other words: it is the number of links originating (or ending, which is the same in an undirected network) in this node. Degree centrality can be seen as a measure of activity (of a node). In this study degree centralities are denoted using the Greek letter delta δ.

## 3. Arrays derived from networks

We recall the definitions of delta, gamma, and alpha-arrays of a network G, as introduced, e.g., by Egghe (2024a,b).

Delta arrays (Egghe, 2024b). The array of degree centralities of a network G with N nodes is denoted as

$$\Delta_G = (\delta_1(G), \delta_2(G), \ldots, \delta_N(G)). \tag{2}$$

where the delta values are ranked decreasingly. Hence, nodes are numbered and named according to their rank in this array. Clearly,

$\sum_{i=1}^{N} \delta_i = 2\,(\#E)$, a notion which is known as the total degree of the network.

Gamma arrays (Egghe, 2024b). We denote the neighboring (or gamma) array of G, as

$$\Gamma(G) = (\gamma_1(G), \gamma_2(G), \ldots, \gamma_N(G)) \ , \tag{3}$$

where the gamma-value of a node i, $\gamma_i(G)$, is equal to the sum of the degree centralities of its zeroth and first-order neighbors, i.e., the degree of the node plus the sum of the degrees of adjacent nodes. Again the elements in array (3) are ordered decreasingly. Consequently, the same index number in (2) or (3) may refer to a different node. The neighboring index of G is defined as $\upsilon(G) = \sum_{i=1}^{N} \gamma_i$(G).

Alpha-arrays (Egghe, 2024a). If $\alpha_j$, j = 1,…, N-1, denotes the number of times distance j occurs in the network G, then the array

$$AF(G) = (\alpha_1(G), \alpha_2(G), \ldots, \alpha_{N-1}(G)) \tag{4}$$

is called the $\alpha$ – array of the network G. Note that for an alpha-array indices do not refer to nodes, but to distances. Consequently, array (4) is not necessarily decreasing.

## 4. The spider family: (M,K,L)

Here we introduce spiders (see Egghe, 2024b): these networks consist of a core of M > 0 nodes, with on each vertex K ≥ 0 “legs” of length L ≥ 0. M, K, and L are natural numbers.

Construction of a spider, denoted as $Sp_{M,K,L}$: a complete network on M nodes acts as the core; then, to each of these M nodes K chains are attached of length L. Informally, we refer to such a set of K chains as a bundle. Hence the number of nodes N = M+MKL, the number of (undirected) links is $\frac{M(M-1)}{2} + MKL$ and the number of shortest paths is

$$\frac{N(N-1)}{2} = \frac{(M+MKL)(M+MKL-1)}{2} = \frac{M(M-1)}{2} - \frac{MKL}{2} + M^2KL + \frac{M^2K^2L^2}{2}\,. \quad (5)$$

If M=1 and L=K=0 then there are no links.

The next figure (Fig.2) shows some spiders.

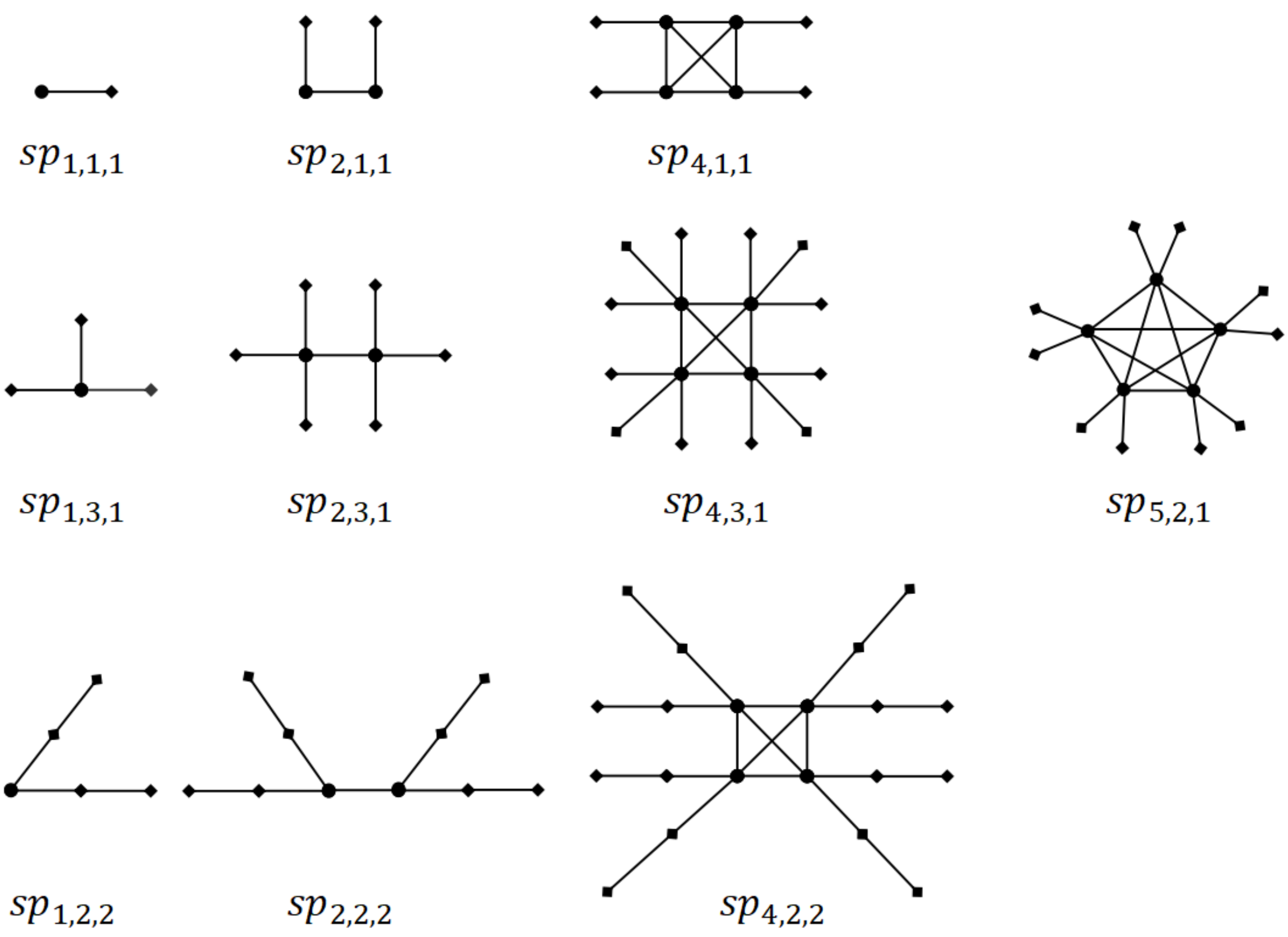

Figure 2: spiders

Special cases

Real spiders (biological spiders have eight legs, usually with seven segments): Hence M = 8, K = 1, L = 7, or M=1, K=8, L=7, depending on how one sees its main body (as a tightly knit construction, or as an abstract “ball”).

If K = 0, or equivalently L = 0, we have a complete network.

If M=1 and L =1 then we have a star (for all K > 0).

If M=1 and L > 1 we have an extended star (for all K > 0).

The case M=1, K = 1, and any L > 0 yields a chain, which is also a rooted tree.

If M=2 and K=1 then we have a chain for every L > 0; every chain with an even number of nodes is of this type.

**5. The arrays of the spider family and other network indicators**

We begin with the calculations of the three arrays (delta, gamma, alpha).

a) The array of degree centralities

The general delta array (M > 1, K > 0, L > 1) is:

$$\Delta_{Sp_{M,K,L}} = \left( \underbrace{M-1+K, \ldots, M-1+K}_{M\ times}, \underbrace{2, \ldots, 2}_{KM(L-1) times}, \underbrace{1, \ldots, 1}_{KM\ times} \right). \quad (6)$$

We check that the total degree is equal to two times the number of links. The total degree is M(M-1+K) + 2KML -2KM + KM = 2KML + M(M-1), which is 2(LKM+M(M-1)/2).

Now we consider some special cases.

If K=0 (equivalently L = 0), then we have a complete network, and only the first group of values in (9) remains (if M>1):

$$\Delta_{Sp_{M,0,L}} = \left( \underbrace{M-1, \ldots, M-1}_{M\ times} \right). \quad (7)$$

If M=1, L=1, and K > 0 then we have a star and the middle part of the array (6) disappears. This leads to:

$$\Delta_{Sp_{1,K,1}} = \left( K, \underbrace{1, \ldots, 1}_{K\ times} \right). \quad (8)$$

If M=1 and L>1 we have an extended star with delta-array:

$$\Delta_{Sp_{1,K,L}} = \left( K, \underbrace{2, \ldots, 2}_{K(L-1) times}, \underbrace{1, \ldots, 1}_{K\ times} \right). \quad (9)$$

If M=1, K=1, L>1 we have a chain

$$\Delta_{Sp_{1,1,L}} = \left(1,\ \underbrace{2,\ldots,2}_{(L-1)times}\ ,1\right).$$

after reranking this becomes: $\Delta_{Sp_{1,1,L}} = \left(\underbrace{2,\ldots,2}_{(L-1)times}\ ,1,1\right).$ (10)

If moreover L =1 we have just one segment with $\Delta_{Sp_{1,1,1}} = (1,1)$

If L =1, K > 0, M > 1, then

$$\Delta_{Sp_{M,K,1}} = \left(\underbrace{M-1+K,\ldots,M-1+K}_{M\ times},\ \underbrace{1,\ldots,1}_{KM\ times}\right). \quad (11)$$

Finally we consider the case M=2, K=1 and L> 0

$$\Delta_{Sp_{2,1,L}} = \left(\underbrace{2,\ldots,2}_{2L\ times}\ ,1,1\right). \quad (12)$$

b) The neighboring array

The general gamma array, i.e. with M > 1, K > 0, L > 3. We consider first the contributions of the core nodes, then those of the leg nodes neighboring the core nodes, then the other leg nodes, except the terminal node one and the one before it, then the one before a terminal node, and finally the terminal leg nodes. This yields:

$$\Gamma_{Sp_{M,K,L}} = \left(\begin{array}{c}\underbrace{M(M-1+K)+2K,\ldots,M(M-1+K)+2K}_{M\ times}, \\ \underbrace{M+K+3,\ldots,M+K+3}_{MK\ times},\ \underbrace{6,\ldots,6}_{(L-3)MK\ times}\ ,\ \underbrace{5,\ldots,5}_{MK\ times}\ ,\ \underbrace{3,\ldots,3}_{MK\ times}\end{array}\right).$$

(13)

Again we consider some special cases.

If K = L = 0, then we have a complete network, and only the first group of values in (13) remains.

$$\Gamma_{Sp_{M,0,L}} = \left(\underbrace{M(M-1),\dots,M(M-1)}_{M\ times}\right). \quad (14)$$

If M=1, L=1, and K > 0 then we have a star, leading to:

$$\Gamma_{Sp_{1,K,1}} = \left(2K, \underbrace{K+1}_{K\ times}\right). \quad (15)$$

If M=1, L = 2 and K>0, then

$$\Gamma_{Sp_{1,K,2}} = \left(3K, \underbrace{K+3}_{K\ times}, \underbrace{3,\dots,3}_{K\ times}\right). \quad (16)$$

If M=1, L = 3 and K > 0, then

$$\Gamma_{Sp_{1,K,3}} = \left(3K, \underbrace{K+4}_{K\ times}, \underbrace{5,\dots,5}_{K\ times}, \underbrace{3,\dots,3}_{K\ times}\right). \quad (17)$$

If M=1, L > 3 and K>0, then

$$\Gamma_{Sp_{1,K,L}} = \left(3K, \underbrace{K+4}_{K\ times}, \underbrace{6,\dots,6}_{K(L-3)\ times}, \underbrace{5,\dots,5}_{K\ times}, \underbrace{3,\dots,3}_{K\ times}\right). \quad (18)$$

If M=1, K=1, L=1, (one segment) then

$$\Gamma_{Sp_{1,1,1}} = (2,2). \quad (19)$$

If M=1, K=1, L = 2 (a 2-chain) then

$$\Gamma_{Sp_{1,1,2}} = (4,3,3). \quad (20)$$

If M=1, K=1, L = 3 (a 3-chain) then

$$\Gamma_{Sp_{1,1,3}} = (5,5,3,3). \quad (21)$$

If M=1, K=1, L >3 (an L- chain), then

$$\Gamma_{Sp_{1,1,L}} = \left(\underbrace{6,\dots 6}_{(L-3)times}, 5,5,3,3\right). \quad (22)$$

If M > 1, K > 0, L=1, then

$$\Gamma_{Sp_{M,K,1}} = \left(\underbrace{M(M-1+K)+K,\dots,M(M-1+K)+K}_{M\ times}, \underbrace{M+K,\dots,M+K}_{MK\ times}\right). \quad (23)$$

If M=2, K=1, L > 0, (a (2L+1)-chain), then

$$\Gamma_{Sp_{2,1,L}} = \left( \underbrace{6, \dots, 6}_{2(L-1)\ times}, 5,5,3,3 \right). \quad (24)$$

Formula (24) also holds for L = 1. We point out that the same chain (with an even number of nodes) has a different gamma array, depending on the fact that there is one core node or there are two core nodes.

c) The frequency of the shortest paths array

The alpha array of a spider is denoted as $AF_{Sp_{M,K,L}} = (\alpha_1, \alpha_2, \dots, \alpha_{N-1})$, where we simply wrote $\alpha_j$ for $\alpha_j(Sp_{M,K,L})$.

It is clear that $\alpha_1 = \frac{M(M-1)}{2} + MKL$. From now on we only consider distances strictly larger than one, hence d(i,j) where i or j does not belong to the core, i.e., the M-complete graph. We further note that if 2L+2 ≤ j ≤ N-1 then $\alpha_j = 0$ (if such values j exist).

We make the following considerations.

(1) There are KM chains in a spider. Within each chain we have once a distance L, twice a distance L-1, three times a distance L-2, and so on ending with L-1 times a distance 2. Recall that distances 1 are already considered above.

(2). Between two chains ending in the same core node, i.e., belonging to the same bundle of K chains, we have in total $M\,\frac{K(K-1)}{2}$ times the distances i+1, …, i+L for all $i \in \{1, \dots, L\}$.

(3). Between two chains of different bundles, i.e.,, ending in different core nodes. We have two cases:

(3a) we have $\frac{MK(KM-K)}{2}$ times the distances i+2, …, i+L+1 for every $i \in \{1, \dots, L\}$, where we did not include the core node for the second chain.

(3b) If we now consider the core node for the second chain we have KM(M-1) times the distances i+1, for all $i \in \{1, \dots, L\}$.

Now we can determine $\alpha_2, \dots, \alpha_{2L+1}$.

$$\alpha_2 = KM(L-1) + M\frac{K(K-1)}{2} + KM(M-1). \qquad (25)$$

Here we have used cases (1), (2), and (3b), where for case (2) only i=1 is possible.

$$\alpha_3 = KM(L-2) + M\frac{K(K-1)}{2}\,2 + KM(M-1) + \frac{MK^2(M-1)}{2}. \qquad (26)$$

Here we used cases (1), (2), (3a), and (3b). Note that for (2) i=2 and i=3 are possible. For (3a) only i=1 is possible.

We can continue in this way for 2 ≤ j ≤ L-1 leading to

$$\alpha_j = KM\big(L-(j-1)\big) + M\frac{K(K-1)}{2}(j-1) + KM(M-1) + \frac{MK^2(M-1)}{2}(j-2). \qquad (27)$$

and ending in the same way with

$$\alpha_L = KM + M\frac{K(K-1)}{2}(L-1) + KM(M-1) + \frac{MK^2(M-1)}{2}(L-2). \qquad (28)$$

Now $$\alpha_{L+1} = M\frac{K(K-1)}{2}\,L + KM(M-1) + \frac{MK^2(M-1)}{2}(L-1). \qquad (29)$$

Here case (1) is not possible anymore. For case (2) i = 1, ..,L are possible, for case (3a) i =1, …, L-1 are possible.

$$\alpha_{L+2} = M\frac{K(K-1)}{2}\,(L-1) + \frac{MK^2(M-1)}{2}\,L\ . \qquad (30)$$

Now case (3b) is not possible anymore. For case (2) i = 2, .., L and for case (3a) i = 1, .., L is possible. Similarly, we have

$$\alpha_{L+3} = M\frac{K(K-1)}{2}\,(L-2) + \frac{MK^2(M-1)}{2}\,(L-1), \qquad (31)$$

and so on, till we reach

$$\alpha_{2L} = M\frac{K(K-1)}{2} + \frac{MK^2(M-1)}{2}\,2. \qquad (32)$$

Here only i = L for case (2) and for case (3a) i =L-1 and i=L are possible (with L ≥ 2).

Finally,

$$\alpha_{2L+1} = \frac{MK^2(M-1)}{2}\,, \qquad (33)$$

where only case (3a) for i=L is possible (with L ≥ 1).

In this way, all alpha values are determined. As this calculation is rather complicated we check if the sum of all alpha-values is equal to the number of different shortest distances, i.e. N(N-1)/2.

$$\sum_{j=1}^{N-1} \alpha_j = \sum_{j=1}^{2L+1} \alpha_j$$

$$= \frac{M(M-1)}{2} + MKL + KM\big((L-1) + \cdots + 1\big)$$

$$+ M\,\frac{K(K-1)}{2}\left(\underbrace{1 + \cdots + (L-1) + L + (L-1) + \cdots + 1}_{= L^2}\right)$$

$$+ KM(M-1)L$$

$$+ M\frac{K^2(M-1)}{2}(\underbrace{1 + \cdots + (L-1) + L + (L-1) + \cdots 1}_{= L^2})$$

$$= \frac{M(M-1)}{2} + MKL + MK\frac{L(L-1)}{2} + M\frac{K(K-1)}{2}\,L^2 + KM(M-1)L + M\frac{K^2(M-1)}{2}L^2$$

Expanding and simplifying yields

$$= \frac{M(M-1)}{2} - \frac{MKL}{2} + M^2KL + \frac{M^2K^2L^2}{2} = \frac{N(N-1)}{2}$$ by (9).

Special cases.

If K = L = 0 then N = M and we have a complete network with

$$\mathrm{AF}_{Sp_{M,0,L}} = \left(\frac{M(M-1)}{2},\ \underbrace{0, \ldots, 0}_{M-2\ times}\right). \qquad (34)$$

If M=1, L=1 and K > 1 then N = K+1, we have a star, leading to:

$$\mathrm{AF}_{Sp_{1,K,1}} = \left(K, \frac{K(K-1)}{2},\ \underbrace{0, \ldots, 0}_{K-1\ times}\right). \qquad (35)$$

If M=1, L = 2 and K>1, then N = 2K+1 and

$$\mathrm{AF}_{Sp_{1,K,2}} = \left(2K, \frac{K(K+1)}{2}, K(K-1),\ \frac{K(K-1)}{2},\ \underbrace{0, \ldots, 0}_{2K-3\ times}\right). \qquad (36)$$

If M=1, K=1, L=1, then N= 2 (one segment) and

$$\mathrm{AF}_{Sp_{1,1,1}} = 1. \quad (37)$$

If M=1, K=1, L = 2 (a 2-chain) then N = 3 and

$$\mathrm{AF}_{Sp_{1,1,2}} = (2,1). \quad (38)$$

If M=1, K=1, L >2 (an L- chain), then N= L+1 and

$$\mathrm{AF}_{Sp_{1,1,L}} = (L, L-1, L-2,\ L-3, \dots, 1\ ). \quad (39)$$

If M=2, K=1, L > 0, (a (2L+1)-chain), then N = 2L+1 and

$$\mathrm{AF}_{Sp_{2,1,L}} = (2L+1, 2L, \dots, 1). \quad (40)$$

We note that formula (44) also holds for L =1.

If M = 2, K=2, L=1, (the “H”), then N= 6, and

$$\mathrm{AF}_{Sp_{2,2,1}} = (5, 6,4,0,0). \quad (41)$$

d) The diameter of a spider

General case (M > 1)

Then the diameter is 2L+1, namely from a terminal node to another terminal node, associated with another core node. This formula is also correct when L = 0.

If M =1 and K > 1, then the diameter is 2L

If M =1 and K =1, then the diameter is L

If M=1 and K=L=0, then the diameter is 0

e) The density of a spider

We already observed that a spider has LKM+M(M-1)/2 edges. As it has M+M*K*L nodes the density of a spider is equal to

$$D(Sp_{M,K,L}) = \frac{2MKL+M(M-1)}{M(1+KL)(M(1+KL)-1\ )} = \frac{2KL+M-1}{M(1+KL)^2-(1+KL)}. \quad (42)$$

Some special cases.

If K=L=0, M > 1, we have a complete network whose density is 1.

If M ( > 1) and K are fixed, then $\lim_{L\to\infty} D(Sp_{M,K,L}) = 0$ and similarly if M (> 1) and L are fixed then $\lim_{K\to\infty} D(Sp_{M,K,L}) = 0$. Both results are intuitively clear.

f) The h-index of a spider

This is by definition (Hirsch, 2005; Schubert et al., 2009) the h-index of the ordered array of degree centralities, hence of the delta-array. Values for the h-index are shown in Table 1. Note that the value M-1 (K=0) also holds for M =1.

Table 1: h-indices of a spider

| Parameters | M > 1 K > 0 | K=0 | M=1, K>0,L=1 | M=1, K > 1, L>1 | M=K=1, L > 2 | M = K =1, 1 ≤ L ≤ 2 |
|---|---|---|---|---|---|---|
| h-index | M | M-1 | 1 | 2 | 2 | 1 |

**6. Dynamical properties: definitions of small worlds**

Two types of small worlds (associated with the alpha and delta arrays), as introduced in earlier work, will be considered (Egghe, 2024a, 2024b; Egghe & Rousseau, 2024). We always assume that a sequence of connected, undirected networks $(G_N)_{N\in\mathbb{N}}$ is given.

6.1 Definitions of small worlds derived from the degree distribution

Because we will define here small worlds derived from the degree distribution we will use the abbreviation DSW (degree small world).

6.1.1 Small worlds derived from the largest degree (DSWL)

If $\Delta_{G_N} = (\delta_1(G_N), \delta_2(G_N), \dots, \delta_N(G_N))$, $N \in \mathbb{N}$, is the delta-array of $\mathrm{G}_N$, then $(\mathrm{G}_N)_{N\in\mathbb{N}}$ is a degree small world based on the largest degree if

$$\lim_{N\to+\infty} \frac{\delta_1(G_N)}{\ln(N)} = +\infty. \quad (43)$$

We informally say that $(\mathrm{G}_N)_{N\in\mathbb{N}}$ is DSWL.

6.1.2 Small worlds derived from the average degree (DSWA)

If $\overline{\delta(G_N)}, N \in \mathbb{N}$, denotes the average degree in network $\mathrm{G}_N$,

$$\overline{\delta(G_N)} = \frac{1}{N} \sum_{j=1}^{N} \delta_j(G_N). \quad (44)$$

then $(\mathrm{G}_N)_{N\in\mathbb{N}}$ is a degree small world based on the average degree if

$$\lim_{N\to+\infty} \frac{\overline{\delta(G_N)}}{\ln(N)} = +\infty. \quad (45)$$

In the same vein as above we say that $(\mathrm{G}_N)_{N\in\mathbb{N}}$ is DSWA.

6.2 Small worlds based on the frequency of distances distribution

6.2.1 Small worlds based on the diameter (SWD)

If $d_N$, $N \in \mathbb{N}$, is the diameter of $\Omega_N$, defined as

$$d_N = max\{d(A,B);\ A, B \in G_N\}, \quad (46)$$

then $(G_N)_{N\in\mathbb{N}}$ is a SWD if there exists a finite constant $C \geq 0$ such that

$$\lim_{N\to+\infty} \frac{d_N}{\ln(N)} = C. \quad (47)$$

Note that $d_N$ is short for diam($\mathrm{G}_N$).

6.2.2 Small worlds based on the average distance (SWA)

If $\mu_N$, $N \in \mathbb{N}$, denotes the average distance between two different elements in $\mathrm{G}_N$:

$$\mu_N = \frac{1}{N(N-1)} \sum_{\substack{A,B \in G_N \\ A \neq B}} d(A,B). \tag{48}$$

then $(G_N)_{N\in\mathbb{N}}$ is an SWA if there exists a finite number $C \geq 0$ such that

$$\lim_{N\to+\infty} \frac{\mu_N}{\ln(N)} = C. \tag{49}$$

**7. Sequences of spiders and their dynamic, small-world properties**

We wonder which sequences of spiders are small worlds? We check the three types of small worlds for M to infinity (K, L fixed), K to infinity (M and L fixed); and L to infinity (K and L fixed).

7.1 DSWLs and spiders

We have to check if

$$\lim_{N\to+\infty} \frac{\delta_1(G_N)}{\ln(N)} = +\infty. \tag{50}$$

where $G_N$ is now a spider with N nodes. Recall that N = M(1+KL). We know that $\delta_1 = M - 1 + K$. Letting M, K, and L tend to infinity, we have to calculate the following limits:

$$\lim_{M\to+\infty} \frac{M-1+K}{\ln(M(1+KL))} = +\infty \,; \lim_{K\to+\infty} \frac{M-1+K}{\ln(M(1+KL))} = +\infty; \lim_{L\to+\infty} \frac{M-1+K}{\ln(M(1+KL))} = 0. \tag{51}$$

This proves that spiders with fixed K and L or with fixed M and K are DSWL. If M and K are fixed, then a spider tends intuitively to an extended star and hence it is not a small world.

7.2 DSWAs and spiders

We have to check if

$$\lim_{N\to+\infty} \frac{\overline{\delta(G_N)}}{\ln(N)} = +\infty. \tag{52}$$

where again $G_N$ is a spider with N = M(1+KL) nodes. Now,

$\overline{\delta(G_N)} = \frac{1}{M(1+KL)}\left(2(\mathrm{LKM} + \frac{\mathrm{M(M-1)}}{2})\right) = \frac{2LKM+M(M-1)}{M(1+KL)} = \frac{M+2KL-1}{1+KL}$, hence we have to calculate:

$\lim_{N\to+\infty} \frac{M+2KL-1}{(1+\mathrm{KL})\ln(M(1+KL))}$ . Again letting M, K, and L to infinity (and keeping the other two variables fixed) leads to:

$\lim_{M\to+\infty} \frac{M+2KL-1}{(1+\mathrm{KL})\ln(M(1+KL))} = +\infty$ ; $\lim_{K\to+\infty} \frac{M+2KL-1}{(1+\mathrm{KL})\ln(M(1+KL))} = 0$ ; and

$\lim_{L\to+\infty} \frac{M+2KL-1}{(1+\mathrm{KL})\ln(M(1+KL))} = 0.$ (53)

This proves that spiders with fixed K and L are DSWA. If M is fixed then spiders are not DSWA.

7.3 SWDs and spiders

We have to check if formula (47), namely

$$\lim_{N\to+\infty} \frac{d_N}{\ln(N)} = C \quad (54)$$

holds, where $d_N$ is short for diam($G_N$). In the general case, the diameter of a spider is 2L+1. Hence we calculate:

$$\lim_{M\to+\infty} \frac{2L+1}{\ln(M(1+KL)} = 0; \lim_{K\to+\infty} \frac{2L+1}{\ln(M(1+KL))} = 0;$$

$$\lim_{N\to+\infty} \frac{2L+1}{\ln(M(1+KL))} = +\infty.. \quad (55)$$

If L is fixed we even have an ultra-small world. Again we do not have a small world for M and K fixed, and L variable.

7.4 SWAs and spiders

Now we have to check if there exists C ≥ 0, such that $\lim_{N\to+\infty} \frac{\mu_N}{\ln(N)} = C$ ,

with $\mu_N = \frac{1}{N(N-1)} \sum_{\substack{A,B \in \mathrm{G}_N \\ A\neq B}} d(A,B)$.

$$\sum_{\substack{A,B \in G_N \\ A \neq B}} d(A,B) = \sum_{j=1}^{N-1} j\alpha_j = \sum_{j=1}^{2L+1} j\alpha_j$$

$$= \frac{M(M-1)}{2} + MKL + KM\left(\underbrace{2(L-1) + 3(L-2) + \cdots + L.1}_{I}\right)$$

$$+ KM(M-1)\left(\underbrace{2 + \cdots + (L+1)}_{II}\right)$$

$$+M\frac{K(K-1)}{2}\left(\underbrace{1.2 + 2.3 + \cdots + L(L+1)}_{III} + \right.$$

$$\left.\underbrace{(L-1)(L+2) + (L-2)(L+3) + \cdots + 1.(2L)}_{IV}\right)$$

$$\frac{K^2(M-1)}{2}\left(\underbrace{1.3 + 2.4 + \cdots + (L-2)L + (L-1)(L+1) + L(L+2)}_{V}\right.$$

$$\left. + \underbrace{(L+3)(L-1) + (L+4)(L-2) \ldots + (2L)(2) + (2L+1).1}_{VI}\right)$$

Now we determine the sums (I) to (VI).

(I) = $\sum_{j=1}^{L-1}(j+1)(L-j) = L\sum_{j=1}^{L-1} j - \sum_{j=1}^{L-1} j^2 + L(L-1) - \sum_{j=1}^{L-1} j$

= $(L-1)\frac{L(L-1)}{2} + L(L-1) - \frac{(L-1)L(2(L-1)+1)}{6}$

$$= \frac{L(L-1)}{6}(3(L-1) + 6 - (2L-1))$$

= $\frac{L(L-1)(L+4)}{6}$;

(II) = $\frac{(L+1)(L+2)}{2} - 1 = \frac{L(L+3)}{2}$

(III) = $\sum_{j=1}^{L} j(j+1) = \sum_{j=1}^{L} j^2 + \sum_{j=1}^{L} j = \frac{L(L+1)(2L+1)}{6} + \frac{L(L+1)}{2}$

= $\frac{L(L+1)}{6}(2L+1+3) = \frac{L(L+1)(L+2)}{3}$

(IV) = $\sum_{j=1}^{L-1}(L-j)(L+j+1) = (L-1)L(L+1) - \sum_{j=1}^{L-1} j^2 - \sum_{j=1}^{L-1} j$

= $(L-1)L(L+1) - \frac{(L-1)L(2L-1)}{6} - \frac{L(L-1)}{2}$

$= \frac{L(L-1)}{6}\,(4L-4) = \frac{2L(L-1)^2}{3}$

(V) = $\sum_{j=1}^{L} j(j+2) = \sum_{j=1}^{L} j^2 + 2\sum_{j=1}^{L} j = \frac{L(L+1)(2L+1)}{6} + 2\frac{L(L+1)}{2}$

$= \frac{L(L+1)}{6}\,(2L+1+6) = \frac{L(L+1)(2L+7)}{6}$

(VI) = $\sum_{j=1}^{L-1}(L+j+2)(L-j) = (L-1)(L+2)L - \sum_{j=1}^{L-1} j^2 - 2\sum_{j=1}^{L-1} j$

$= (L-1)(L+2)L - \frac{(L-1)L(2L-1)}{6} - 2\,\frac{L(L-1)}{2}$

$= \frac{L(L-1)}{6}\,(6(L+2)-(2L-1)-6) = \frac{L(L-1)(4L+7)}{6}$

Hence we have:

$$\sum_{j=1}^{2L+1} j\alpha_j = \frac{M(M-1)}{2} + MKL + KM\frac{L(L-1)(L+4)}{6} + KM(M-1)\frac{L(L+3)}{2}$$

$$+ M\frac{K(K-1)}{2}\left(\frac{L(L+1)(L+2)}{3} + \frac{2L(L-1)^2}{3}\right) + M\frac{K^2(M-1)}{2}\left(\frac{L(L+1)(2L+7)}{6} + \frac{L(L-1)(4L+7)}{6}\right)$$

This is a closed form. Dividing by N(N-1) = $M^2(1+KL)^2$ yields the average $\mu_N$ . Then

$$\lim_{M\to+\infty}\frac{\mu_N}{\ln(N)} = 0;\ \lim_{K\to+\infty}\frac{\mu_N}{\ln(N)} = 0;\ \lim_{L\to+\infty}\frac{\mu_N}{\ln(N)} = +\infty;$$

This shows that we also have an ultra-small world in the SWA sense except for the case of growing L-values and fixed M and K (leading to a network that looks like an extended star).

## 8. Conclusion

We have shown that spiders provide examples of small worlds. We suggest that spiders can act as models for core-periphery studies (Borgatti and Everett, 1999) and as such wonder if real core-periphery networks too, behave as small worlds. Besides as a model for core-periphery networks, spiders are also interesting examples in more formal network studies as they provide a fluent transition between two extreme situations, namely chains and complete networks.